# On Number Representation

Rafael I. ROFA

June, 2013

Perth, Western Australia

Email: rrofa12@bigpond.com

## ABSTRACT

Place value numbers, such as the binary or decimal numbers can be represented by the end vertices (leaf or pendant vertices) of rooted symmetrical trees. Numbers that consist of at most a fixed number of digits are represented by vertices that are equidistant from the root vertex and the corresponding number representations do not depend on the distance from the root vertex. In this paper, we introduce place value number systems which are representable by rooted symmetrical trees and in which the representation of a number depends on the distance of the corresponding vertex from the root vertex. Such dependence activates the role of zero in such a way as to render its function equivalent to that of any other single digit number. Thus, in addition to being a place value holder, the digit zero (just as any other single digit numeral) affects the value of a number regardless of its position. For example 012 is different, in the new systems, from 12. As such, these new number systems could be thought of as a natural development for the role of zero. We also illustrate how addition is performed in these newly constructed number systems. In addition to being mathematical structures which could be of mathematical interest, these new number systems could possibly have applications in computing and computing security.

## 1. Introduction

Number systems, particularly the decimal and binary systems are key structures upon which most of mathematics and science rests.

The base b place value number system, where b is a positive integer, consists of numbers of the form

$$(a_n a_{n-1} \ldots a_0)_b = \sum_{k=0}^{n} a_k b^k$$

For example, in the base 10 (decimal) number system, we have

21738=(21738)$_{10}$=2x10$^4$+1x10$^3$+7x10$^2$+3x10$^1$+8x10$^0$

Historically, place value number systems developed throughout several civilizations. For a comprehensive history of these developments refer to [1], [2], [3], [4], and [5]. It is known that the Babylonians were among the first to use a place value number system as far back as the 19$^{th}$ century BC [1]. They used a sexagesimal (base 60) system. However, the Babylonian system lacked the numeral 0 (at a later stage Babylonians used two slanted wedges to keep an empty place value) [3]. The numeral 0, both as a place value holder and as a number by itself appeared in the Hindu-Arabic number system [6]. The Hindu-Arabic number system is the decimal place value number system that is still in use today together with similar various other systems having different bases such as the binary number system which has 2 as its base instead of the base 10 decimal. However, the role of zero in these systems is not equivalent to that of other nonzero numerals. For example, in the decimal number system, 012 and 12 have the same value; the zero in 012 has no effect on the value of the number. In general, a zero that is positioned to the left of any nonzero numeral does not affect the value of the number that it appears in. This property, in these systems, renders the role of the numeral zero not equivalent to that of any other nonzero numeral. In this paper, we develop new place value number systems in which the role of the numeral zero is equivalent to that of any other nonzero numeral.

Consider the set of all binary numbers each of which consisting of at most k digits. These numbers can be enumerated using a rooted symmetrical tree with k levels and the $2^k$ such binary numbers are represented by the end vertices (leaf vertices) of such a tree. For k=3, figure 1 shows the 8 binary numbers being represented by the 8 end vertices (labeled 000, 001, 010, 011, 100, 101, 110, and 111) of a rooted symmetrical tree with root vertex v.

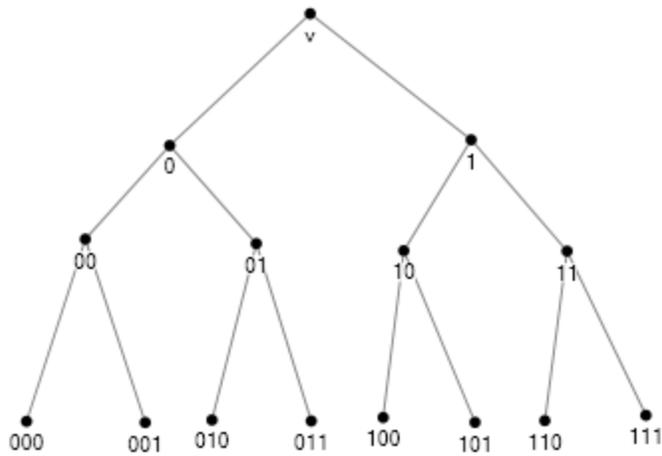

Created with NodeXL (http://nodexl.codeplex.com)

**Figure 1**

There are 8 binary numbers each consisting of at most 3 digits (notice that 001 for example is equal to 1). In comparison there are 14 vertices, excluding root vertex, in the tree that is used to enumerate the 8 binary numbers. The place value number systems developed in this paper, utilize all the vertices of a rooted symmetrical tree in such a way that any two vertices represent *distinct* numbers. For example, there will be 14 different "binary" numbers represented by the tree of figure 1, instead of the known 8 binary numbers.

## 2. Definitions and Preliminary Results

**Definition 2.1** Let $K = (k_1, k_2, ..., k_r)$ be a sequence of positive integers. For all q in $\{2,...,r,r+1\}$, let $b_q = 1 + k_{q-1} + k_{q-1}k_{q-2} + ... + k_{q-1}k_{q-2}... k_1$. The set $B = \{b_q \,/\, 1 \leq q \leq r\}$, where $b_1 = 1$, is said to be the place value set generated by K and K is said to be a base or radix sequence.

**Lemma 2.2** Let $B = \{b_1, b_2, ..., b_r\}$ be the place value set generated by the radix sequence $K = (k_1, k_2, ..., k_r)$, then

For all q in $\{2,...,r+1\}$, we have

$$b_q = 1 + k_{q-1}b_{q-1}$$

**Proof.**

By Definition 2.1, we have, $\forall q \in \{2,...r+1\}$

$b_q = 1 + k_{q-1} + k_{q-1}k_{q-2} + ... + k_{q-1}k_{q-2}...k_1$

$\quad = 1 + k_{q-1}(1 + k_{q-2} + k_{q-2}K_{q-3} + ... + k_{q-2}k_{q-3}...k_1)$

$\quad = 1 + k_{q-1}b_{q-1}$

**Example 2.3** Let $K = (4,3,2)$, then by Definition 2.1

$b_1 = 1$

$b_2 = 1 + k_1 = 1 + 4 = 5$

$b_3=1+k_2+k_2k_1 =1+3+12=16$

$b_4=1+k_3+k_3k_2+k_3k_2k_1=1+2+6+24=33$

Hence, the place value set B generated By K is {1,5,16}

Alternatively, by Lemma 2.2(ii), we have

$b_2=1+k_1b_1=1+4=5$

$b_3=1+k_2b_2=1+3\times 5=16$, and

$b_4=1+k_3b_3=1+2\times 16=33$

The following Lemma is a generalization of Lemma 2.2 and will be useful in Section 3.

**Lemma 2.4** Let $B=\{b_1,b_1,…,b_r\}$ be the place value set generated by the radix sequence $K=(k_1,k_2,…,k_r)$. For all integers q in {2,…,r+1} and all integers s in {1,…,q}, we have

$$b_q = (k_{q-1}-1)b_{q-1}+(k_{q-2}-1)b_{q-2}+…+(k_{q-s+1}-1)b_{q-s+1}+k_{q-s}b_{q-s}+s$$

**Proof**

We prove Lemma 2.4 by the Principle of Mathematical Induction on q. If q=2, then s=1, and Lemma 2.4 asserts that

$b_2=k_1b_1+1$, which is, by Lemma 2.2, is true.

Let p be an integer such that 2≤p≤r+1, and assume that Lemma 2.4 is true for q=p. In other words, assume that for all s in {0,1,…,p-1}, we have

$b_p=(k_{p-1}-1)b_{p-1}+(k_{p-2}-1)b_{p-2}+…+(k_{p-s+1}-1)b_{p-s+1}+k_{p-s}b_{p-s}+s$ (1)

Now, by Lemma 2.2 part, we have

$b_{p+1}=1+k_p b_p$

$\quad =1+(k_p-1)b_p+b_p$ \hfill (2)

Substituting the expression for $b_p$ from equation (1) into equation (2) we get

$b_{p+1}=1+(k_p-1)b_p+(k_{p-1}-1)b_{p-1}+...+(k_{r-s-1}-1)b_{s+1}+k_{r-s}b_s+p-s$

$\quad =1+(k_p-1)b_p+(k_{p-1}-1)b_{p-1}+...+(k_{p-s+1}-1)b_{p-s+1}+k_{p-s}b_{p-s}+s$

$\quad =(k_p-1)b_p+(k_{p-1}-1)b_{p-1}+...+(k_{p-s+1}-1)b_{p-s+1}+k_{p-s}b_{p-s}+s+1$

Thus, if Lemma 2.4 is true for q=p, where 2≤p, then it is also true for q=p+1, Hence, By The Principle of Mathematical Induction Lemma 2.4 is true for all q in {2,…,r+1}. This completes the proof.

Notice that if s=1 then q-s=q-1, hence for all q in {2,…,r+1} we have $b_q=k_{q-1}b_{q-1}+1$, which is what Lemma 2.2 states. Thus, Lemma 2.2 is a special case of Lemma 2.4. Furthermore, if q-s=1, then $b_{s-q}=b_1=1$, hence

$b_q = (k_{q-1}-1)b_{q-1}+(k_{q-2}-1)b_{q-2}+...+(k_1-1)+q$

**Example 2.6** Let K=(10,10,10,10,10), then $b_1=1$, $b_2=11$ $b_3=111$ $b_4=1,111$ $b_5=11,111$ and $b_6=111,111$.

We apply Lemma 2.4 for the following cases:

(a) For q=5 and s=3, we have
$b_5=(k_4-1)b_4+(k_3-1)b_3+k_2 b_2+3$
$\quad =9\times 1,111+9\times 111+10\times 11+3=11,111$

(b) For q=4 and s=1, we have
$b_4 = k_3 b_3 + 1$
$= 10 \times 111 + 1 = 1{,}111$

(c) For q=5 and s=4, we have
$b_5 = (k_4 - 1)b_4 + (k_3 - 1)b_3 + (k_2 - 1)b_2 + k_1 b_1 + 4$
$= 9 \times 1{,}111 + 9 \times 111 + 9 \times 11 + 10 \times 1 + 4 = 11{,}111$

In the next section, we present the main result of this paper.

## 3. Distance sensitive Number Representation

Let $B = \{b_1, b_2, \ldots, b_r\}$ be the place value set generated by the radix sequence $K = (k_1, k_2, \ldots, k_r)$. Let R be $\{(i_r, i_{r-1}, \ldots, i_{r-q}) / q \in \{0, 1, \ldots, r-1\}$ and $\forall x \in \{1, \ldots, r\}, 0 \leq i_x \leq k_x - 1\}$. Note that the elements of R correspond to the set of all vertices of a rooted symmetrical tree T with degree sequence $k_r, k_{r-1}, \ldots k_1$. An element, say $(i_r, \ldots, i_{r-q})$, of R represents a vertex at distance $q+1$ from the root vertex. Using counting techniques, we have

$|R| = k_r + k_r k_{r-1} + \ldots + k_r k_{r-1} \ldots k_1 = b_{r+1} - 1$.

Next, we state and prove the main theorem of this paper.

**Theorem 3.1** Let K, B and R be as defined above. The function F defined by $F((i_r, i_{r-1} \ldots, i_{r-q})) = i_r b_r + i_{r-1} b_{r-1} \ldots + i_{r-q} b_{r-q} + q$, for $(i_r, \ldots, i_{r-q}) \in R$, is injective. If $F((i_r, \ldots, i_{r-p})) = n$, where $(i_r, \ldots, i_{r-p}) \in R$, then the sequence $(i_r, \ldots, i_{r-p})$ is said to *be the Distance Sensitive Number Representation of n with respect to the place value set B, and is denoted by $R_B(n)$. The set R is said to be the Distance Sensitive number system associated with the place value set B. Furthermore, p is said to be the distance term of $R_B(n)$*

**Proof** Let M be the range of F. We prove Theorem 3.1 by showing that the inverse relation is a function. So let $m \in M$, then there exists $(i_r,...,i_{r-q})$, for some q in $\{0,...,r-1\}$ satisfying the conditions of Theorem 3.1 such that $m=F((i_r,...,i_{r-q}))= i_r b_r+...+i_{r-q}b_{r-q}+q$. The algorithm below, shows that $(i_r,...,i_{r-q})$ is unique and therefore, it shows that the inverse of F is a well defined function.

### The Distance Sensitive Division Algorithm.

let $m \in M=$Range(F), then as mentioned above, there exists a sequence $(i_r,...,i_{r-q})$, for some q in $\{0,1,...,r-1\}$ such that $F((i_r,...,i_{r-q}))=i_r b_r+...+i_{r-q}b_{r-q}+q=m$ the following algorithm yields the unique sequence $R_B(m)=(i_r,...,i_{r-q})$.

- Evaluate $m/b_r$:

$m/b_r=(i_r b_r+...+i_{r-q}b_{r-q}+q)/b_r$

$=i_r$, remainder $i_{r-1}b_{r-1}+...+i_{r-q}b_{r-q}+q<b_r$ [By Lemma 2.4 and $i_x \leq k_x$ for $x \in \{1,...,r\}$].

Set $n_1= i_{r-1}b_{r-1}+...+i_{r-q}b_{r-q}+q$

If $n_1=0$, set $R_B(m)=(i_r)$. Algorithm terminates.

If $n_1 \neq 0$, set $m_1=n_1-1$

- Evaluate $m_1/b_{r-1}$:

$m_1/b_{r-1}=i_{r-1}$, remainder $i_{r-2}b_{r-2}+...+i_{r-q}b_{r-q}+q-1<b_{r-1}$ [By Lemma 2.4]

Set $n_2 = i_{r-2}b_{r-2} + \ldots + i_{r-q}b_{r-q} + q - 1$

If $n_2 = 0$, set $R_B(m) = (i_r, i_{r-1})$. Algorithm terminates.

If $n_2 \neq 0$, set $m_2 = n_2 - 1$

.
.
.

- Evaluate $m_{r-q+1}/b_{r-q+1}$ :

$M_{r-q+1}/b_{r-q+1} = i_{r-q+1}$ remainder $i_{r-q}b_{r-q} + 1 < b_{r-q+1}$ [By Lemma 2.4]

Set $n_q = i_{r-q}b_{r-q} + 1$

If $n_q = 0$, Set $R_B(m) = (i_r, i_{r-1}, \ldots, i_{r-q+1})$. Algorithm terminates.

If $n_q \neq 0$, set $m_q = n_q - 1$ $(= i_{r-q}b_{r-q})$

- Evaluate $m_q/b_{r-q}$:

$m_q/b_{r-q} = i_{r-q}$, remainder 0.

Set $R_B(m) = (i_r, i_{r-1}, \ldots, i_{r-q})$

Algorithm terminates.

Thus the Algorithm terminates and for each $m \in \text{Range}(F)$ it yields a unique value $R_B(m)$ such that $F^{-1}(m) = R_B(m)$. Hence, $F^{-1}$ is well defined and it follows that F is injective.

This completes the proof of Theorem 3.1

Note: Where there is no confusion, the sequence representing $R_B(m)$ (that is ($i_r, i_{r-1}, \ldots, i_{r-q}$)) is sometimes written as $i_r i_{r-1} \ldots i_{r-q}$, otherwise it is written as $(i_r i_{r-1} \ldots i_{r-q})_B$, and brackets could be omitted as well. This last notation for $R_B(m)$ differentiates it from other number representations such as the decimal or binary

representations and indicates its relevance to the particular place value set B.

**Lemma 3.2** Let M be the range of the function F, then

i)    $|M|=b_{r+1}-1$, and
ii)    $M=\{0,1,...,b_{r+1}-2\}$

Proof.i) By Theorem 3.2, F: R→M is injective, hence

(i) $|M|=|R|=k_r+k_rk_{r-1}+...+k_rk_{r-1}...k_1=b_{r+1}-1$.

ii) The minimal element of M occurs when q=0 and $i_r=0$ which is 0.

Moreover, the maximal element n occurs when q=r-1 and $i_r, i_{r-1},...,i_1$ have their maximum values, hence

$n=(k_r-1)b_r+(k_{r-1}-1)b_{r-1}+...(k_1-1)b_1+r-1=b_{r+1}-2$ [By Lemma 2.5]

But by (i) $|M|=b_{r+1}-1$, hence,

$M=\{0,1,...,b_{r+1}-2\}$.

This completes the proof of Lemma 3.2

The function F as defined in Theorem 3.1 and the Distance Sensitive Division Algorithm are used to convert Distance Sensitive representation of numbers into decimals and vice versa, respectively. The following is an example.

**Example 3.3** Let K=(4,5,2,6) be a radix sequence. Find

(a)    The place value set B generated by K.
(b)    Range(F)
(c)    $F((5)_B)$

(d) $F((1,0,3)_B)$

(e) $F((1,1,4)_B)$

(f) $F((2,0,1)_B)$

(g) $R_B(265)$

(h) $R_B(70)$

(i) $R_B(101)$

(j) $R_B(113)$

Note: If it is not otherwise indicated, all numbers used are ternary or decimal numbers.

**Solution**:

(a) $b_1=1$

$b_2=k_1b_1+1=5$

$b_3=k_2b_2+1=26$

$b_4=k_3b_3+1=53$, and

$b_5=k_4b_4+1=319$

Hence, $B=\{53,26,5,1\}$

(b) By Lemma 3.2 we have Range(F)=$\{0,1,\ldots,317\}$

(c) $F((5)_B)=5b_4=265$

(d) $F((1,0,3)_B)=b_4+0b_3+3b_2+2=70$

(e) $F((1,1,4)_B)=b_4+b_3+4b_2+2=101$

(f) $F((2,0,1)_B)=2b_4+0b_3+b_2+2=113$

(g) Applying The Distance Sensitive Division Algorithm:

(265)/$b_4=i_4=5$. Remainder $n_1=0$. Hence, $R_B(265)=5_B$

Alternatively, $R_B(265)=F^{-1}(265)=5_B$ [Using part c]

(h) Applying The Distance Sensitive Division Algorithm:

(70)/$b_4=i_4=1$, remainder $n_1=17$.

Set $m_1=n_1-1=16$.

$m_1/b_3=16/26=i_3=0$, remainder $n_2=16$.

    Set $m_2=n_2-1=15$

    $m_2/b_2=15/5=i_2=3$, remainder $n_3=0$. Hence, $R_B(70)=103_B$

    Alternatively, $R_B(70)= F^{-1}(70)=103_B$ [By part d]

(i)   Applying The Distance Sensitive Division Algorithm:

    $(101)/b_4=101/53=i_4=1$, remainder $n_1=48$.

    Set $m_1=n_1-1=47$.

    $m_1/b_3=47/26=i_3=1$, remainder $n_2=21$.

    Set $m_2=n_2-1=20$

    $m_2/b_2=20/5=i_2=4$, remainder $n_3=0$. Hence, $R_B(101)=114_B$

(j)   $R_B(113)=F^{-1}(113)=201_B$ [By part f]

Or alternatively, The Distance Sensitive Division Algorithm can be used to find $R_B(113)$.

**Definition 3.4** A distance sensitive decimal number system (or a distance decimal number system) is a distance sensitive number system with radix sequence whose entries are all tens.

**Note 3.5** Definition 3.4 can be extended to include all other standard place value number systems. For example, a distance binary number system is a distance sensitive number system with radix sequence whose entries are all 2's

**Example 3.6** Let $R_d$ be the distance decimal number system associated with the radix sequence D of size d (By Definition 3.3, all elements of D are tens). Then the elements of the place value set B generated by D are: $b_1=1$, $b_2=11$, $b_3=111$, … $b_d=11…1$ (1 repeated d times). Distance decimal numbers are represented with reference to d, for example, if $d=2$, $(23)_D$ is written as $23_2$. Table 1 below shows the decimal numbers from 0 to 109 together with the corresponding distance decimal representations or numbers.

| Decimal | D.D | Decimal | D.D | Decimal | D.D | Decimal | D.D |
|---|---|---|---|---|---|---|---|
| 0 | 0 | 28 | 25 | 56 | 50 | 84 | 76 |
| 1 | 00 | 29 | 26 | 57 | 51 | 85 | 77 |
| 2 | 01 | 30 | 27 | 58 | 52 | 86 | 78 |
| 3 | 02 | 31 | 28 | 59 | 53 | 87 | 79 |
| 4 | 03 | 32 | 29 | 60 | 54 | 88 | 8 |
| 5 | 04 | 33 | 3 | 61 | 55 | 89 | 80 |
| 6 | 05 | 34 | 30 | 62 | 56 | 90 | 81 |
| 7 | 06 | 35 | 31 | 63 | 57 | 91 | 82 |
| 8 | 07 | 36 | 32 | 64 | 58 | 92 | 83 |
| 9 | 08 | 37 | 33 | 65 | 59 | 93 | 84 |
| 10 | 09 | 38 | 34 | 66 | 6 | 94 | 85 |
| 11 | 1 | 39 | 35 | 67 | 60 | 95 | 86 |
| 12 | 10 | 40 | 36 | 68 | 61 | 96 | 87 |
| 13 | 11 | 41 | 37 | 69 | 62 | 97 | 88 |
| 14 | 12 | 42 | 38 | 70 | 63 | 98 | 89 |
| 15 | 13 | 43 | 39 | 71 | 64 | 99 | 9 |
| 16 | 14 | 44 | 4 | 72 | 65 | (100) | 90 |
| 17 | 15 | 45 | 40 | 73 | 66 | (101) | 91 |
| 18 | 16 | 46 | 41 | 74 | 67 | (102) | 92 |
| 19 | 17 | 47 | 42 | 75 | 68 | (103) | 93 |
| 20 | 18 | 48 | 43 | 76 | 69 | (104) | 94 |
| 21 | 19 | 49 | 44 | 77 | 7 | (105) | 95 |
| 22 | 2 | 50 | 45 | 78 | 70 | (106) | 96 |
| 23 | 20 | 51 | 46 | 79 | 71 | (107) | 97 |
| 24 | 21 | 52 | 47 | 80 | 72 | (108) | 98 |
| 25 | 22 | 53 | 48 | 81 | 73 | (109) | 99 |
| 26 | 23 | 54 | 49 | 82 | 74 | | |
| 27 | 24 | 55 | 5 | 83 | 75 | | |

**Table 1**

Notice that the decimal numbers corresponding to the distance decimals ranging from 90 to 99 are 3 digit decimal numbers.

**Example 3.7** Let R be the distance binary number system with radix sequence K=(2,2,2,). The place value set generated by K is B={7,3,1}. Table 2, below, displays all of the eight binary numbers each consisting of at most 3 digits together with the corresponding distance binary numbers whose radix set is K. (Again notice that there are 14 distance binary numbers each having at most 3 digits compared to only 8 binary numbers each consisting of at most 3 digits).

| Decimal | Binary | Distance Binary |
| --- | --- | --- |
| 0 | 0 | 0 |
| 1 | 1 | 00 |
| 2 | 10 | 000 |
| 3 | 11 | 001 |
| 4 | 100 | 01 |
| 5 | 101 | 010 |
| 6 | 110 | 011 |
| 7 | 111 | 1 |
| 8 | - | 10 |
| 9 | - | 100 |
| 10 | - | 101 |
| 11 | - | 11 |
| 12 | - | 110 |
| 13 | - | 111 |

Table2

Notice that the lengths of the distance binary numbers [that is the number of digits in each number] follow repetitive cycles related to $b_1$, $b_2$, and $b_3$. Moreover, the number of distinct distance binary numbers is $2^3+2^2+2+1$ which is the number of vertices (excluding the root vertex) of the rooted symmetrical tree in figure 1 [Lemma 3.2]

A number N, however large, can be expressed as a decimal by first finding the largest power of 10 that goes into N and then applying the usual division algorithm. In a similar manner, it is possible to express a decimal number as a distant decimal (or a distant binary, etc.) number by using the Distance Division Algorithm which is stated above. The following example illustrates the process.

**Example 3.8** Express the decimal number N=205,556 as a distance decimal number.

**Solution** First notice that the place value set of a distance decimal number system is B={1,11,111,…}, where the number of elements of B is equal to the number of entries (each of which is a 10) in the radix sequence that generates B. Hence, to find $R_B(205,556)$, one has to find the largest element of B that is contained in 205,556 which is 111,111. Hence we set $b_r=b_6=111,111$ and therefore, B={111111,11111,1111,111,11,1}. Applying the Distance Division Algorithm we find that $R_B(205,556)=(185)_B$

In this example, Notice that B is the smallest place value set in which it is possible to express the decimal 201,503 as a distance decimal number and as such we call it *the critical place value set with respect to*

*N*. On the other hand, Any distance decimal number system with a place value set B' that contains B will contain $R_{B'}(205,556)$. For example, if $B_1=B\cup\{1111111\}$, then, by the Distance Division Algorithm, $R_{B1}(205,556)=(0184999)_{B1}$. If $B_2=B_1\cup\{11111111\}$, then $R_{B2}(205,556)=(00184998)_{B2}$. Table 3 summarizes these results and extends them to progressively larger sets, $B_3, B_4,...$

| $B_i$ | $R_{Bi}(205,556)$ |
|---|---|
| B | $(185)_B$ |
| $B_1$ | $(0184999)_{B1}$ |
| $B_2$ | $(00184998)_{B2}$ |
| $B_3$ | $(000184997)_{B3}$ |
| $B_4$ | $(0000184996)_{B4}$ |
| $B_5$ | $(00000184995)_{B5}$ |
| $B_6$ | $(000000184994)_{B6}$ |
| $B_7$ | $(0000000184993)_{B7}$ |
| $B_8$ | $(00000000184992)_{B8}$ |
| $B_9$ | $(000000000184991)_{B9}$ |
| $B_{10}$ | $(0000000000184990)_{B10}$ |
| $B_{11}$ | $(0000000000018499)_{B11}$ |
| $B_{12}$ | $(00000000000184989)_{B12}$ |
| $B_{13}$ | $(000000000000184988)_{B13}$ |
| $B_{14}$ | $(0000000000000184987)_{B14}$ |
| $B_{15}$ | $(00000000000000184986)_{B15}$ |
| $B_{16}$ | $(000000000000000184985)_{B16}$ |
| $B_{17}$ | $(0000000000000000184984)_{B17}$ |
| $B_{18}$ | $(00000000000000000184983)_{B18}$ |
| $B_{19}$ | $(000000000000000000184982)_{B19}$ |
| $B_{20}$ | $(0000000000000000000184981)_{B20}$ |
| $B_{21}$ | $(00000000000000000000184980)_{B21}$ |
| $B_{22}$ | $(000000000000000000000018498)_{B22}$ |

Table 3

## 4. Addition of distance sensitive numbers

In this section we illustrate the process of finding the sum of two distance decimal numbers. The place value set of a distance decimal number system is B={1,11,111,…} [See example 3.6]. Furthermore, each digit of a distance decimal number is between 0 and 9 inclusive. Adhering to the notation of example 3.6, a distance sensitive number is given with reference to the size of the associated radix sequence. For example the distance sensitive number $23_2$ is equivalent to the decimal 26 (2x11+3x1+1) whereas $23_4$ is equivalent to the decimal 2556 (2x1111+3x111+1). In the case when the number of digits of a distance decimal number is equal to size of the radix sequence of the distance decimal set to which it belongs, the suffix is completely omitted; for example, $328_3$ is written as 328 and $0251_4$ is written as 0251.

To add two distance decimal numbers we perform the following steps.

1) <u>Step 1</u>: Addition is performed by adding corresponding equal place value digits of each number. Moreover, the "distance" terms [the distance term in each number is one less than the number of the rest of the digits of that number] are added. For example, to add $32_3$ and 15 we notice that the place value of the digit 3 in the first number is $b_3=111$, the place value of the digits 2 and 1 in the first and second numbers, respectively, is $b_1=11$, and the place value of the digit 5 in the second number is $b_1=1$. Moreover, the "distance" term in both numbers is equal to one. Hence,

$32_3$ =  3  2    1
11   =       1  1  1
Sum =   3  3  1  2

Thus,

$32_3+11=331$.

The decimal equivalent of the above equation is $356+13=369$

2) <u>Step 2</u>: Adding two equal place value single digits as in step 1 above, could result in a double digit number, hence a "carrying" procedure is required. Such a procedure is done in accordance with Lemma 2.2. For example

32   =3   2   1
28   =2   8   1
Sum =5  (10)  2
      =6    0   1   [By Lemma 2.2 (ii) $10b_1 +1=b_2$]

Hence,

$32+28=60$.

In decimal numbers, this is equivalent to $36+31=67$

3) <u>Step 3</u>: Performing steps 1 and 2 above could result in a sequence of digits in which the last digit (the distance digit) is not one less than the number of the rest of the digits. This violates Theorem 3.1. In such cases a number is "borrowed" or "lended" to get the required value of the distance term. The following examples illustrate these processes.

**Example 3.1**

```
231   =2  3   1   2
52    =      5   2   1
Sum   =2  8   3   3
      =2  8   4   2   [1 is subtracted from the distance
                       Term, making it conform with the
                       distance rule, and added to the
                       coefficient of $b_1$=1]
```

Thus,
231+52=284.
Or in decimals,
258+58=316.

**Example 3.2**

```
$25_5$  = 2   5                   1
324   =           3   2   4   2
Sum   = 2   5   3   2   4   3
      = 2   5   3   2   3   4   [1 is taken from the coefficient
                                 of $b_1$ and is added to the
                                 "distance" term to make it
                                 conform with the distance
                                 rule.]
```

Thus,
$25_5$ + 324 = 25323.
In decimal numbers, the above equation is equivalent to

27778 + 361 = 28139

**Example 3.3**

| | | | | |
|---|---|---|---|---|
| $3_4$ = | 3 | | | 0 |
| $26_3$ = | | 2 | 6 | 1 |
| Sum = | 3 | 2 | 6 | 1 |
| = | 3 | 2 | 5 | (11) 1 [one $b_2$=11 is borrowed and placed as a $b_1$ coefficient] |
| = | 3 | 2 | 5 | 9 3 [2 are borrowed from the $b_1$ coefficient and added to distance term to make it conform to the distance rule] |

Thus,

$3_4 + 26_3 = 3259$.

The decimal equivalent to the above equation is

3333 + 289 = 3622

**Example 3.4**

| | | | | |
|---|---|---|---|---|
| 099 = | 0 | 9 | 9 | 2 |
| 1 = | | | 1 | 0 |
| Sum = | 0 | 9 | (10) | 2 |
| = | 0 | (10) | 0 | 1 [$10b_1+1=b_2$] |
| = | 1 | 0 | 0 | 0 [$10b_2+1=b_3$] |
| = | 1 | | | 0 [the zero coefficients of $b_2$ And $b_1$ are retracted making the distance term conform the distance rule] |

Thus,

099 + 1 = $1_3$

The decimal equivalence of the above equation is

110 + 1 = 111

**Example 3.5**

| | | | | | |
|---|---|---|---|---|---|
| 3529 = | 3 | 5 | 2 | 9 | 3 |
| 11 = | | | 1 | 1 | 1 |
| Sum = | 3 | 5 | 3 | (10) | 4 |
| = | 3 | 5 | 4 | 0 | 3 |

$[10b_1+1=b_2]$. Notice that the zero coefficient of $b_1$ must not be retracted as in the previous example because the distance term is already in agreement with the distance rule.

Example 4.1

Let R be the distance decimal number system with K=(10,10). Find

(a) $(3)_R + (5)_R$

(b) $(23)_R + (15)_R$

(c) $(21)_R + (7)_R$

(d) $(12)_R + (15)_R$

(e) $(27)_R+(33)_R$

(f) $(25)_R+(37)_R$

(g) $(5)_R+(6)_R$

Solution

The place value set generated by K in $\{11,1\}$ with $b_1=11$ and $b_2=1$. Addition is done mod $b_0-1$ which is 110.

(a) $(3)_R+(5)_R=(8)_R$ [Decimal equivalent: 33+55=88]

The detailed procedure is as follows:

$(3)_R=3b_1$

$(5)_R=5b_1$

Sum=$8b_1$ [Rule 1]

(b) $(23)_R+(15)_R=(39)_R$. [Decimal equivalent: 26+17=43]

The detailed procedure is as follows:

$(23)_R=2b_1+3b_2+1$

$(15)_R=b_1+5b_2+1$

Sum=$3b_1+8b_2+2$ [Rule 1]

$=3b_1+9b_2+1$ [Rule 3, constant is 1 less than the total number of the other terms, so we 1 from the constant 2 and add it to the Coefficient of $b_2$, remembering that $b_2=1$].

$=(39)_R$

(c) $(21)_R+(7)_R=(91)_R$ [Decimal equivalent:77+24=101]

　　Detailed procedure:

　　$(21)_R=2b_1+b_2+1$

　　$(7)_R=7b_1$

　　Sum=$9b_1+b_2+1$ [Rule 1]

　　　　=$(91)_R$

(d) 　$(12)_R+(15)_R=(28)_R$ [Decimal equivalent:14+17=31]

　　Detailed procedure:

　　$(12)_R=b_1+2b_2+1$

　　$(15)_R=b_1+5b2+1$

　　Sum=$2b_1+7b_2+2$ [Does not satisfy Rule 3]

　　　　=$2b_1+8b_2+1$ [Satisfies Rule 3]

　　　　=$(28)_R$

(e) $(27)_R+(33)_R=(60)_R$ [Decimal equivalent:30+37=67]

　　Detailed procedure:

　　$(27)_R=2b_1+7b_2+1$

　　$(33)_R=3b_1+3b_2+1$

　　Sum=$5b_1+10b_2+2$ [Coefficient of $b_2$ doesn't satisfy Rule 4]

　　　　=$6b_1+0b_2+1$　　[Rule 2: By lemma 2.2, $k_2b_2+1=b_1$, or

　　　　　　　　　　　　$10b_2+1=b_1$. Notice that 1 has been

　　　　　　　　　　　　deducted from the constant term

　　　　　　　　　　　　and added to $10b_2$ to make one $b_1$]

　　　　=$(60)_R$

(f) $(25)_R+(37)_R=(62)_R$ [Decimal equivalent:28+41=69]

　　Detailed procedure:

　　$(25)_R=2b_1+5b_2+1$

　　$(37)_R=3b_1+7b_2+1$

Sum=$5b_1+12b_2+2$ [Rules 3 and 4 not satisfied]

=$6b_1+2b_2+1$   [By Lemma 2.2, $10b_2+1=b_1$]

=$(62)_R$

(g) $(5)_R+(7)_R=(2)_R$ [Decimal equivalent:55+77=22(mod110)

Detailed procedure:

$(5)_R=5b_1$

$(7)_R=7b_1$

Sum=$12b_1$ [does not satisfy Rule 4]

=$2b_1$   [Applying Rule 5, addition is done mod $b_0-1$ which, by Lemma 2.2 is $10b_1$]

**References**

1. G Ifrah, *A universal history of numbers: From prehistory to the invention of the computer* (London, 1998).
2. R Kaplan, *The nothing that is: a natural history of zero* (London, 1999).
3. L C Karpinski, *The history of arithmetic* (New York, 1965).
4. K W Menninger, *Number words and number symbols: A cultural history of numbers* (Boston, 1969).
5. D E Smith and L C Karpinski, *The Hindu-Arabic numerals* (Boston, 1911).